\documentstyle[amssymb,12pt]{amsart}

\newtheorem{prop}{Proposition}[section]
\newtheorem{prop:def}{Proposition-Definition}[section]

\newtheorem{lemma}{Lemma}[section]
\newtheorem{thm}{Theorem}[section]
\newtheorem{cor}{Corollary}[section]

\theoremstyle{remark}
\newtheorem{remark}{Remark}

\begin{document}

\newcommand{\nc}{\newcommand} \nc{\on}{\operatorname}

\nc{\¦}{{|}}

\nc{\pa}{\partial}

\nc{\cA}{{\cal A}}\nc{\cB}{{\cal B}}\nc{\cC}{{\cal C}}
\nc{\cE}{{\cal E}} \nc{\cG}{{\cal G}}\nc{\cH}{{\cal
H}}\nc{\cX}{{\cal X}}\nc{\cR}{{\cal R}} \nc{\cL}{{\cal
L}}\nc{\cK}{{\cal K}}\nc{\cO}{{\cal O}} \nc{\cF}{{\cal
F}}\nc{\cM}{{\cal M}} \nc{\cW}{{\cal W}}\nc{\cV}{{\cal V}}

\nc{\sh}{\on{sh}}\nc{\Id}{\on{Id}}\nc{\Diff}{\on{Diff}}
\nc{\ad}{\on{ad}}\nc{\Der}{\on{Der}}\nc{\End}{\on{End}}
\nc{\res}{\on{res}}\nc{\ddiv}{\on{div}} \nc{\FS}{\on{FS}}
\nc{\card}{\on{card}}\nc{\dimm}{\on{dim}}\nc{\gr}{\on{gr}}
\nc{\Jac}{\on{Jac}}\nc{\Ker}{\on{Ker}} \nc{\Den}{\on{Den}}
\nc{\Imm}{\on{Im}}\nc{\limm}{\on{lim}}\nc{\Ad}{\on{Ad}}
\nc{\ev}{\on{ev}} \nc{\Hol}{\on{Hol}}\nc{\Det}{\on{Det}}
\nc{\Cone}{\on{Cone}} \nc{\pseudo}{{\on{pseudo}}}
\nc{\class}{{\on{class}}}\nc{\rat}{{\on{rat}}}
\nc{\local}{{\on{local}}}\nc{\an}{{\on{an}}}
\nc{\Lift}{{\on{Lift}}}\nc{\Mer}{{\on{Mer}}}\nc{\mer}{{\on{mer}}}
\nc{\lift}{{\on{lift}}}\nc{\diff}{{\on{diff}}}\nc{\Aut}{{\on{Aut}}}
\nc{\DO}{{\on{DO}}}\nc{\Frac}{{\on{Frac}}}\nc{\cl}{{\on{class}}}
\nc{\Fil}{{\on{Fil}}}

\nc{\Bun}{\on{Bun}}\nc{\diag}{\on{diag}}\nc{\KZ}{{\on{KZ}}}
\nc{\CB}{{\on{CB}}}\nc{\out}{{\on{out}}}\nc{\Hom}{{\on{Hom}}}
\nc{\FO}{{\on{FO}}}

\nc{\al}{\alpha}\nc{\de}{\delta}\nc{\si}{\sigma}\nc{\ve}{\varepsilon}\nc{\z}{\zeta}
\nc{\vp}{\varphi}
\nc{\la}{{\lambda}}\nc{\g}{\gamma}\nc{\eps}{\epsilon}
\nc{\PsiDO}{\Psi\on{DO}}\nc{\om}{\omega}

\nc{\AAA}{{\mathbb A}}\nc{\CC}{{\mathbb C}}\nc{\NN}{{\mathbb N}}
\nc{\PP}{{\mathbb P}}\nc{\RR}{{\mathbb R}}\nc{\VV}{{\mathbb V}}
\nc{\ZZ}{{\mathbb Z}}

\nc{\bla}{{\mathbf \lambda}} \nc{\bv}{{\mathbf v}}
\nc{\bz}{{\mathbf z}}\nc{\bt}{{\mathbf t}}
\nc{\bP}{{\mathbf P}} \nc{\kk}{{\mathbf k}}

\nc{\A}{{\mathfrak a}}\nc{\B}{{\mathfrak b}}\nc{\G}{{\mathfrak g}}
\nc{\HH}{{\mathfrak h}}\nc{\mm}{{\mathfrak m}}\nc{\N}{{\mathfrak n}}
\nc{\SG}{{\mathfrak S}}\nc{\La}{\Lambda}

\nc{\wt}{\widetilde}
\nc{\wh}{\widehat} \nc{\bn}{\begin{equation}}\nc{\en}{\end{equation}}
\nc{\SL}{{\mathfrak{sl}}}\nc{\ttt}{{\mathfrak{t}}}
\nc{\GL}{{\mathfrak{gl}}}

\title[Integrable systems associated with elliptic algebras]
{Integrable systems associated with elliptic algebras}

\begin{abstract}
We construct some new Integrable Systems (IS) both classical and
quantum associated with elliptic algebras. Our constructions are
partly based on the algebraic integrability mechanism given by the
existence of commuting families in skew fields and partly - on the
internal properties of the elliptic algebras and their
representations. We give some examples to make an evidence how
these IS are related to previously studied.
\end{abstract}

\author{A. Odesskii and V. Rubtsov}

\address{A.O.: Laboratoire de Physique Th\'eorique, UMR 8627 du CNRS, Universit\'e Paris XI,
91405, Orsay, France. Permanent address: Landau Institute of
Theoretical Physics, 117234, Kosygina, 2 Moscow, Russia}
\address{V.R.: D\'epartement de Math\'ematiques, UMR 6093 du CNRS, Universit\'e d'Angers,
49045 Angers, France and ITEP, 25, Bol. Tcheremushkinskaya,
117259, Moscow, Russia}

\maketitle

\subsection*{Introduction}

This paper is an attempt to establish a direct connection between
two close subjects of modern Mathematical Physics - the theory of
Integrable Systems (IS) and the Elliptic Algebras. The aim of this
connection is two-fold: we clarify some our recent results in the
both domains and fill the natural gap proving that some large
class of the Elliptic Algebras carries in fact the structure of an
IS.

 We will start with a short account in the subject of the
story and will describe briefly a type of the IS's under
cosideration.

 In \cite{Skew}, B.Enriquez and second author had proposed a
construction of commuting families of elements in skew fields.
They explained how to use this construction in Poisson fraction
fields to give an another proof of the integrability of Beauville
- Mukai integrable systems associated with a K3 surface $S$
(\cite{Beau1}).

The Beauville-Mukai systems had appeared as the Lagrangian
fibrations which have the form $S^{[g]} \to \¦\cL \¦ =
\PP(H^0(S,\cL))$, where $S^{[g]}$ is the Hilbert scheme of $g$
points of $S$, equipped with a symplectic structure introduced in
\cite{Mukai}, and $\cL$ is a line bundle on $S$.  Later, the
authors of \cite{Donagi} explained that these systems are natural
deformations of the "separated" (in the sense of \cite{GNR})
versions of Hitchin's integrable systems, more precisely, of their
description in terms of spectral curves (already present in
\cite{Hitchin}). Beauville-Mukai systems can be generalized to
surfaces with a Poisson structure (see \cite{Bottacin}). When $S$
is the canonical cone $\Cone(C)$ of an algebraic curve $C$ these
systems coincide with the separated version of Hitchin's systems.
The paper \cite{Skew} shows how the commuting families
construction provides a quantization of these systems on the
canonical cone.

A quantization of Hitchin's system was proposed in \cite{BD}. It
seems interesting to construct quantization of Beauville-Mukai
systems and to compare it with Beilinson - Drinfeld quantization.
We had conjectured the correspondence between the results of
\cite{BD} and a quantization of fraction fields in \cite{ER}. A
part of this program was realized in the \cite{ER} for the case of
$S = T^*{{\PP}^{1}_{n}},$ where ${\PP}^{1}_{n}$ means a rational
curve with $n$ marked points.

Another main object of the paper is the families of Elliptic
Algebras. These algebras (with 4 generators) were appeared in the
works of Sklyanin \cite{Skl1},\cite{Skl2} and then they were
generalized (for any number of generators) and intensively studied
by Feigin and one of the
authors(\cite{FO1},\cite{FO2},\cite{FO3}). These algebras can be
considered as a flat deformation of the polynomial rings such that
the linear ( in the deformation parameter ) term is given by the
quadratic Poisson brackets. The
geometric meaning of the Poisson structures was underscored in
\cite{FO4}(see also \cite{Pol} and \cite{Mar})  - it is the
Poisson structures of moduli spaces of holomorphic bundles on the
elliptic curve. We will use the recent survey \cite{Od1} as a main
source of the results and references in the theory of elliptic
algebras.

The relevance of Elliptic Algebras to the theory of IS was
manifested from the very beginning. They appeared in Sklyanin's
demonstration of the integrability of the Landau-Lifshitz model
(\cite{Skl1},\cite{Skl2})in the frame of the Faddeev's school
ideology (Quantum Inverse Scattering Method and $R-$ matrix
approach). Later Cherednik had observed the relation of the
Elliptic Algebras defined in \cite{FO1} with the Belavin
$R-$matrix (see \cite{Cher}). An interesting observation of
Krichever and Zabrodin giving an interpretation of a generator in
the Sklyanin's Elliptic Algebra as a Hamiltonian of 2-point
Ruijsenaars IS (\cite{KZ}) was generalized later in \cite{BGOR} to
the case of double-elliptic 2-point classical model.  However, all
applications of these algebras to the IS theory had an indirect
character until the last two years.

We should mention also the results which were obtained in the paper of
Sokolov-Tsyganov (\cite{SokTsyg}) where they construct,(using the
Sklyanin's definition of quadratic Poisson structures), some
classical commuting families associated with this Poisson
algebras. The integrability of these families is implied by the
Sklyanin's ideology of Separated Variables (SoV) and technically
is based on some generalization of the classical methods going
back to Jacobi, Liouville and St$\ddot{a}$ckel ( which 
ideologically is very close to the classical part of the theorems in 
\cite{Skew} and \cite{BabTal} ). However, all their results are in the
"non-elliptic" case.

 Recently, an example of integrable systems associated with
 linear and quadratic Poisson brackets given by the elliptic Belavin-Drinfeld
 classical $r$-matrix was proposed in \cite{KhLO}. This system (an
 elliptic rotator) appears both in finite and infinite-dimensional
 cases. They give an elliptic version of $2D$ ideal
 hydrodynamics on the symplectomorphism group of the 2-dimensional
 torus as well as on a non-commutative torus. We should also
 mention another appearance of the elliptic algebras in the
 context of Non-Commutative geometry (see \cite{CDV}). It would be
 interesting to relate them to a numerous modern attempts to
 define a Non-Commutative version of IS theory.

 In this paper we construct some IS which appear directly in the
 frame of the elliptic algebras.  The elliptic algebras are
 figured here in two-fold way: sometimes, they are carrying the
 commuting families of Hamiltonians, sometimes, - they provide a
 necessary background to our constructions which use their
 properties and representation theory ( basically, the so-called
 "functional realization" and the "bosonization" mappings).

 A family of coordinated (compatible) elliptic Poisson structures
was introduced in \cite{Od2}. This family contains three quadratic
Poisson brackets such that their generic linear combination is the
quasi-classical limit $q_{n}(\cE)$ of an elliptic algebra
$Q_n(\cE,\eta)$. The famous Lenard - Magri scheme provides the
existence of a classical integrable system associated with the
elliptic curve but it was not clear how to get a
quantum counterpart of this system because of lack of the
knowledge how to quantize in general the Magri-Lenard scheme.
Nevertheless, this system can be quantized for $n=2m$ using the approach
developing in \cite{Skew} and it is one of the main results of our
paper.

Some of the "elliptic" commuting elements which we construct in this
article are related to the quantum version of the above-mentioned
bi-hamiltonian system.
Some other families are associated with a special choice
of the elliptic algebra. These families are obtained, in one hand,
as the direct application of the construction from \cite{Skew} to
the elliptic algebras and, in other hand, by using the properties
of the "bosonization" homomorphism, constructed in earlier works
(\cite{FO2},\cite{FO3}). Some of these families (under the
appropriate choice of their numeric parameters) may be interpreted
as examples of algebraic completely integrable systems. We give a
geometric interpretation to some of them describing a link with
the Lagrangian fibrations on symmetric products of elliptic curve
cone, giving a version of the Beauville-Mukai systems (see
\cite{Beau1},\cite{GNR},\cite{Skew},\cite{Mark}).

Roughly speaking, the integrable systems of the first type have as
the phase space a $2m-$ dimensional component of the moduli space of 
parabolic rank two bundles on the given elliptic curve
$\cE$. More precisely, the coordinate ring of the open dense part
of the latter has a structure of a quadratic Poisson algebra
isomorphic to $q_{2m}(\cE)$. We have explicitly verified that the
quantum commuting elements from our construction are the same as
the latter obtained from the Lenard-Magri scheme for $m=3$ (the
first non-trivial case).

Our main theorem (Thm.3.1) takes place for $n=2m$, but we are sure
that there are some interesting integrable quantum systems in the
case of $n=2m+1$. 
It would be interesting to study the
bi-hamiltonian structures giving the algebra $q_{2m+1}(\cE)$ using
the results of
Gelfand-Zakharevich (\cite{GZ}) on the geometry of bi-hamiltonian
systems in the case of impair-dimensional Poisson manifolds. 
The precise quantum version of these systems in the
context of the elliptic algebras $Q_{2m+1}(\cE,\eta)$ is still
obscure and should be a subject of further studies.

The theorems from \cite{Skew} may be also interpreted as an
algebraic version of the SoV method ( as it was argued in
\cite{BabTal}). Hence, it is very plausible that some of our
quantum commuting families arising from the generalization of the
Jacobi-Liouville-St$\ddot{a}$ckel conditions ( which are
guaranteed by existence conditions of the Cartier - Foata NC
determinants) are the quantum elliptic versions of the IS from
\cite{SokTsyg}. We hope to return to this question in our future
paper.

We give also some low-dimensional examples of our construction.

\section{Commuting families in some non-commutative algebras}

Let $A$ be an associative algebra with unit. We will suppose further that
we work with a skew field.

Fix a natural number $n\geq 2$. We will suppose that there are $n$
subalgebras $B_n \subset A$ in $A$ such that for any couple of the
indices $i\neq j, 1\leq i,j \leq n$ the elements $b_{(i)}\in B_i$
and $b_{(j)}\in B_j$ are commuting (and inside of the each
subalgebra $B_i$ the elements generically are not commuting).

Let us consider the following data: take an $n\times (n+1)$ matrix
$\cM$

\begin{equation}\label{var:matrix}
\left(
\begin{array}{ccccc}
b_{(1)}^{0}& b_{(1)}^{1}& \ldots& b_{(1)}^{n-1}& b_{(1)}^{n}\\
b_{(2)}^{0}& b_{(2)}^{1}& \ldots& b_{(2)}^{n-1}& b_{(2)}^{n}\\
\vdots & \vdots& \vdots &   \vdots &         \vdots \\
b_{(n)}^{0}& b_{(n)}^{1}& \ldots& b_{(n)}^{n-1}& b_{(n)}^{n}
\end{array}
\right)
\end{equation}

such that all elements of $i-$th row belong to the $i-$th
subalgebra $B_{(i)}$.

We will denote by $\cM^{i}$ a $n\times n$ matrix which is resulted
from the matrix $\cM$ with the $i-th$ row erased. The
corresponding determinants (if they exist) will be denote by
$M^i$.

 Now we observe that all principal $n+1$ minors $$ M^{0}, M^{1},
\ldots, M^{n} $$
 of the highest order $n$ are correctly defined with the
help of \emph{Cartier - Foata non-commutative determinant}.

It's definition repeats \emph{verbatim} the standard one because
in all $n$ matrices $$ {\cM}^{0},{\cM}^{1}, \ldots, {\cM}^{n} $$
the entries in the \emph{different rows} are commute and each
summand in the standard definition of the determinant contains the
product of $n$  elements of different rows whose product is order
- independent.

The following theorem was proved in \cite{Skew}
\begin{thm}\label{thm:comm}
Suppose that one of the matrix ${\cM}^{i}, i=1,\ldots,n$, say,
${\cM}^0$, is invertible. Then the elements $H_i =
(M^{0})^{-1}M^{i},i=1,\ldots,n$ are commute in the skew field $A$
(which we are identifying with its image in the fraction field
$\Frac(A)$ under the monomorphic embedding $A\to \Frac(A)$.
\end{thm}

The proof of the theorem is achieved by some tedious but a
straightforward induction procedure.

 The similar results were obtained in the framework of so-called multi-parametric spectral
problems in Operator Analysis (\cite{Leeds}) and in the framework
of so called Seiberg-Witten integrable systems associated with a
hyperelliptic spectral curves in \cite{BabTal}.

The important step in the demonstration of the \ref{thm:comm} is
the following \emph{"triangle"} relations which are similar to the
usual Yang -Baxter relation:

\begin{equation}\label{YB1}
M^{i}(M^{0})^{-1}M^{j}=M^{j}(M^{0})^{-1}M^{i}
\end{equation}
\begin{equation}\label{YB2}
B^{ij}(M^{0})^{-1}B^{kj}=B^{ik}(M^{0})^{-1}B^{ij},
\end{equation}
where $B^{ij}$ is the co-factor of the matrix element $b^{i}_{j},\
0\leq i,j \leq n .$

 The theorem \ref{thm:comm}
can be re-formulated to give the following result

\begin{cor}\label{comm:fij}
Let $A$ be an algebra, $(f_{i,j})_{0\leq i \leq n, 1\leq j \leq
n}$ be elements of $A$ such that $$ f_{i,j}f_{k,\ell} = f_{k,\ell}
f_{i,j} $$ for any $i,j,k,\ell$ such that $j\neq\ell$. For any $I
\subset \{0,\ldots,n\},J\subset \{1,\ldots,n\}$ of the same
cardinality, we set $\Delta_{I,J} = \sum_{\sigma\in \on{Bij}(I,J)}
\eps(\sigma) f_{i,\sigma(i)}$. Here $\on{Bij}(I,J)$ denotes a set
of bijections between two sets of the indices $I$ and $J$.
 Assume that the $\Delta_{I,J}$ are all invertible. Set $\Delta_i =
\Delta_{\{1,\ldots,n\},\{0,\ldots,\check i, \ldots,n\}}$ to be the
element corresponded to $J$ with the index $i$ omitted. Then the
$$ H_i = (\Delta_0)^{-1} \Delta_i $$ all commute together.
\end{cor}

\subsection{Poisson commuting families}

We will fix a base field $\kk$ of characteristic $\neq 2$. The
following observation is straightforward:
\begin{lemma} \label{prol:poisson}
If $B$ is an integral Poisson algebra, then there is a
unique Poisson structure on $\Frac(B)$ extending the Poisson
structure of $B$.
\end{lemma}

This structure is uniquely defined by the relations

$$
\{1/f,g\}= - \{f,g\} / f^2, \{1/f,1/g\} =\{f,g\} / (f^2 g^2)
$$.

\medskip

Theorem \ref{thm:comm} has a Poisson counterpart.

\begin{thm} \label{thm:Poisson}
Let $A$ be a Poisson algebra. Assume that $A$ is integral, and let
$\pi : A \hookrightarrow \Frac(A)$ be its injection in its
fraction field. For each $n-$uple  of Poisson subalgebras
$B_1,\ldots,B_n$ of $A$ such that the elements of pair-wise
different subalgebras $B_i$ are Poisson commuting (for any pair of
indices $i,j, i\neq j$ the elements $b_{i}\in B_i$ and $b_{j}\in B_j$ we
have  $\{b_{i}, b_{j}\}=0$.) We will write the analogue of the
matrix (\ref{var:matrix})  like the vector-row: $\cM =
[b^{0},b^{1},\ldots, b^{n}]$ , where $$
 b^{i} =
\left(
\begin{array}{c}
b_{1}^{i}\\ b_{2}^{i}\\
\vdots\\
 b_{n}^{i}
\end{array}
\right).
 $$
 We set $$ \Delta_i^\cl =
\det[b^0,\ldots,\check b^i, \ldots, b^n]. $$  Here as usual we
denote by $\check b^i$ the $i-$th omitted column. Then if
$\Delta_0^\cl$ is nonzero we set $H_i^\cl = \Delta_i^\cl /
\Delta_0^\cl$ and the family $(H_i^\cl)_{i = 1,\ldots,n}$ is
Poisson-commutative: $$ \{H_i^\cl,H_j^\cl\} = 0 $$ for any pair
$(i,j)$.
\end{thm}

\begin{remark}
The elements $b_{i}^{k}$ and $b_{j}^{l}$ of the matrix $\cM$
belongs to the different subalgebras $B_i$ and $B_j$ if $i\neq j$
and hence are Poisson-commute. This condition reminds the
classical constrains on the Poisson brackets between  matrix
elements appeared in XIX century in the papers of St$\ddot
{a}$ckel on the Separation of Variables of Hamilton-Jacobi systems
(\cite{Stak}). So our theorem can be considered as an algebraic
re-definition of the St$\ddot{a}$ckel conditions
\end{remark}

\subsubsection{Plucker relations}
We want to remind the important step of the second proof in
\cite{Skew} which shows the relations between the commuting
elements and the Plucker-like equations.

We have to prove
\begin{equation} \label{var:jac}
\Delta_i^\cl \{\Delta_j^\cl,\Delta_k^\cl\} + \on{cyclic\
permutation\ in\ }(i,j,k) = 0.
\end{equation}
We have $$ \Delta_i^\cl = \sum_{p=1}^n \sum_{\al = 0}^n
(-1)^{p+\al} (b^{\al})^{(p)} ({\Delta}_{\al,i})^{(1\ldots\check p
\ldots n)} , $$ where (if $\al \neq i$)

 $$
  {\Delta}_{\al,i}^{(1\ldots\check p\ldots n)} = (-1)^{1_{\al<i}}
\det[b^{0}\ldots \check b^{\al} \ldots \check b^{i} \ldots
b^{n}]^{(p)} $$

(which means that the $p-$th row in the matrix $[b^0,\ldots,\check
b^{\alpha},\ldots,\check b^{i}\ldots b^n ]$ should be erased.)

We set $1_{\al<i} = 1$ if $\al<i$ and $0$ otherwise. If $\al = i$
we assume ${\Delta}_{\al,i} = 0$. Now we have $$
\{\Delta_i^\cl,\Delta_j^\cl\} = \sum_{p=1}^n \sum_{\al,\beta =
0}^n (-1)^{\al+\beta} (\{b^{\al},b^{\beta}\})^{(p)}
\left(
\Delta_{\al,i}\Delta_{\beta,j} - \Delta_{\beta,i}\Delta_{\al,j}
\right)^{(1\ldots \check p \ldots n)}, $$

so identity (\ref{var:jac}) is a consequence of
\begin{equation} \label{aux}
\forall (i,j,k,\al,\beta,\gamma), \; \sum_{\sigma\in
\on{Perm}(i,j,k)} \eps(\sigma) \Delta_{\al,\sigma(i)}
\Delta_{\beta,\sigma(j)} \Delta_{\gamma,\sigma(k)} = 0.
\end{equation}

When $\card\{\al,\ldots,k\} = 3$, this identity follows from the
antisymmetry relation $\Delta_{i,j} + \Delta_{j,i} = 0$.

When $\card\{\al,\ldots,k\} = 4$ (resp., $5,6$), it follows from
the following Plucker identities (to get (\ref{aux}), one should
set $V = (A^{\otimes n})^{\oplus n}$ and $\La$ some partial
determinant).

Let $V$ be a vector space. Then

-- if $\Lambda \in \wedge^2(V)$, and $a,b,c,d\in V$, then

\begin{equation} \label{pluck1}
 \La(a,b) \La(c,d) - \La(a,c)\La(b,d) + \La(a,d) \La(b,c) = 0;
\end{equation}

-- if $\La\in \wedge^3(V)$ and $a,b,c,b',c' \in V$, then

\begin{align*}
& \La(b,c,c')\La(a,c,b')\La(b,b',c')
+ \La(b,c,b')\La(c,b',c')\La(a,b,c')\\
& - \La(b,c,b')\La(a,c,c')\La(b,b',c')
- \La(b,c,c')\La(c,b',c')\La(a,b,b')
= 0;
\end{align*}

-- if $\La\in \wedge^4(V)$ and $a,b,c,a',b',c'\in V$, then

\begin{align} \label{toprove}
& \nonumber \La(b,c,b',c')\La(a,c,a',c')\La(a,b,a',b')
+ \La(b,c,a',c')\La(a,c,a',b')\La(a,b,b',c') \\
& \nonumber
+ \La(b,c,a',b')\La(a,c,b',c')\La(a,b,a',c') \\
& \nonumber
- \La(b,c,b',c')\La(a,c,a',b')\La(a,b,a',c')
- \La(b,c,a',b')\La(a,c,a',c')\La(a,b,b',c') \\
& - \La(b,c,a',c')\La(a,c,b',c')\La(a,b,a',b') = 0.
\end{align}

We refer to \cite{Skew} for the proof of these identities. We will
need them below in some special situation arising with the
commuting elements in associative and Poisson algebras which are
directly connected with elliptic curves and vector bundles on
them. This Plucker relations can  be interpreted as a kind of
Riemann-Fay identities which are in its turn related to integrable
(difference) equations in Hirota bilinear form.

\medskip

\section{Elliptic algebras}

Now we should describe one of the main heroes of our story - the
elliptic algebras. We will follow to the survey \cite{Od1} in our
notations and also we will refer to it as a main source of the
results and its proofs in this section.

\subsection{Definition and the main properties}

The elliptic algebras are the associative quadratic algebras
$Q_{n,k}(\cE,\eta)$ which were introduced in the papers
\cite{FO1,FO3}. Here $\cE$ is an elliptic curve and $n,k$ are
integer numbers without common divisors ,such that $1\leq k < n$
while $\eta$ is a complex number and $Q_{n,k}(\cE,0) =
\CC[x_1,...,x_n]$.

Let $\cE = \CC/\Gamma$ be an elliptic curve defined by a lattice
$\Gamma = \ZZ\oplus\tau\ZZ, \tau \in {\CC}, \Im \tau > 0$. The
algebra $Q_{n,k}(\cE,\eta)$ has generators $x_i, i\in
{\ZZ}/n{\ZZ}$ subjected to the relations
 $$
\sum_{r\in{\ZZ}/n{\ZZ}}{\frac{\theta_{j-i+r(k-1)}(0)}{\theta_{j-i-r}(-\eta)\theta_{kr}(\eta)}}x_{j-r}x_{i+r}
= 0
$$
and have the following properties: 1) $Q_{n,k}(\cE,\eta) =
{\CC}\oplus Q_1\oplus Q_2\oplus ...$ such that
$Q_{\alpha}*Q_{\beta}=Q_{\alpha + \beta}$, here $*$ denotes the
algebra multiplication. In other words, the algebras
$Q_{n,k}(\cE,\eta)$ are $\ZZ$ - graded;

2) The Hilbert function of $Q_{n,k}(\cE,\eta)$ is $\sum_{\alpha \geq
  0}\dim Q_{\alpha}t^{\alpha} = \frac{1}{(1 - t)^n}$.

  We consider here the theta-functions $\theta_{i}(z), i = 1,\ldots, n$
  as a base in the space of the theta-functions $\Theta_{n}(\Gamma)$ of the order $n$
   which are subordinated to the following
relations of quasi-periodicity $$ \theta_{i}(z+1) = \theta_{i}(z),
\, \theta_{i}(z + \tau) = (-1)^n
\exp(-2\pi\sqrt(-1)nz)\theta_{i}(z), i=0,\ldots,n-1. $$

The theta-function of the order 1 $\theta(z) \in \Theta_{1}(\Gamma)$
satisfies to the conditions $\theta(0)=0$ and $\theta(-z) =
\theta(z + \tau) = -\exp(-2\pi\sqrt(-1)z)\theta(z)$.

We see that the algebra $Q_{n,k}(\cE,\eta)$ for fixed $\cE$ is a
flat deformation of the polynomial ring $\CC[x_1,...,x_n]$. The
linear (in $\eta$)term of this deformation gives rise to a
quadratic Poisson algebra $q_{n,k}(\cE)$.

The geometric meaning of the algebras $Q_{n,k}$ was underscored in
\cite{FO4},\cite{Pol} where it was shown that the quadratic
Poisson structure in the algebras $q_{n,k}(\cE)$ associated with
the above-mentioned deformation is nothing but the Poisson
structure on ${\PP}^{n-1} = {\PP}Ext^1(E,{\cO})$, where $E$ is a
stable vector bundle of rank $k$ and degree $n$ on the elliptic
curve $\cE$.

In what follows we will denote the algebras $Q_{n,1}(\cE,\eta)$ by
$Q_{n}(\cE,\eta)$

\subsection{Algebra $Q_{n}(\cE,\eta)$}

\subsubsection{Construction}
For any $n\in\NN$, any elliptic curve $\cE=\CC/\Gamma$, and any
point $\eta\in\cE$ we construct a graded associative algebra
$Q_n(\cE,\eta)=\CC\oplus F_1\oplus  F_2\oplus\dots$, where
$F_1=\Theta_{n}(\Gamma)$ and
$F_\alpha=S^\alpha\Theta_{n}(\Gamma)$. By
construction, $\dim
F_\alpha=\frac{n(n+1)\dots(n+\alpha-1)}{\alpha!}$. It is clear
that the space $F_\alpha$ can be realized as the space of
holomorphic symmetric functions of $\alpha$ variables
$\{f(z_1,\dots,z_\alpha)\}$ such that
\begin{equation}
\begin{aligned}
f(z_1+1,z_2,\dots,z_\alpha)&=f(z_1,\dots,z_\alpha),\\
f(z_1+\tau,z_2,\dots,z_\alpha)&=(-1)^ne^{-2\pi
inz_1}f(z_1,\dots,z_\alpha).
\end{aligned}
\end{equation}

For $f\in  F_\alpha$ and $g\in  F_\beta$ we define the symmetric
function $f*g$ of $\alpha+\beta$ variables by the formula

$$
 f*g(z_1,\dots,z_{\alpha+\beta})=\frac{1}{\alpha!\beta!}
\sum_{\sigma\in S_{\alpha+\beta}}
f(z_{\sigma_1} + \beta\eta,\ldots,z_{\sigma_\alpha} + \beta\eta)
g(z_{\sigma_{\alpha+1}}-\alpha\eta,\ldots,z_{\sigma_{\alpha+\beta}}-
\alpha\eta)\times
$$
\begin{equation}
\times{\prod\limits_{\tiny
\begin{array}{c}
1\leq i\leq\alpha\\
\alpha+1\leq j\leq\alpha+\beta
\end{array}
}}
 \frac{\theta(z_{\sigma_i}-z_{\sigma_j}-n\eta)}
{\theta(z_{\sigma_i}-z_{\sigma_j})}.
\end{equation}

 In particular, for $f,g\in F_1$ we have $$
f*g(z_1,z_2)=f(z_1 + \eta)g(z_2-\eta)\frac{\theta(z_1-z_2-n\eta)}
{\theta(z_1-z_2)}+f(z_2 + \eta)g(z_1-\eta)\frac{\theta(z_2-z_1-n\eta)}
{\theta(z_2-z_1)}. $$ Here $\theta(z)$ is a theta function of
order one.

\begin{prop}
If $f\in  F_\alpha$ and  $g\in F_\beta$, then $f*g\in
F_{\alpha+\beta}$. The operation $*$ defines an associative
multiplication on the space $\oplus_{\alpha\ge0}F_\alpha$
\end{prop}

\subsubsection{Main properties of the algebra $Q_n(\cE,\eta)$}

By construction, the dimensions of the graded components of the
algebra $Q_n(\cE,\eta)$ coincide with those for the polynomial
ring in $n$ variables. For $\eta=0$ the formula for $f*g$ becomes
$$
f*g(z_1,\dots,z_{\alpha+1})=\frac1{\alpha!\beta!}\sum_{\sigma\in
S_{\alpha+\beta}}f(z_{\sigma_1},\dots,z_{\sigma_\alpha})
g(z_{\sigma_{\alpha+1}},\dots,z_{\sigma_{\alpha+\beta}}). $$ This
is the formula for the ordinary product in the algebra
$S^*\Theta_{n}(\Gamma)$, that is, in the polynomial ring in $n$
variables. Therefore, for a fixed elliptic curve $\cE$ (that is,
for a fixed modular parameter $\tau$) the family of algebras
$Q_n(\cE,\eta)$ is a deformation of the polynomial ring. In
particular, there is a Poisson algebra, which we denote by
$q_n(\cE)$. One can readily obtain the formula for the Poisson
bracket on the polynomial ring from the formula for $f*g$ by
expanding the difference $f*g-g*f$ in the Taylor series with
respect to $\eta$. It follows from the semi-continuity arguments
that the algebra $Q_n(\cE,\eta)$ with generic $\eta$ is determined
by $n$ generators and $\frac{n(n-1)}2$ quadratic relations. One
can prove (see \S2.6 in \cite{Od1}) that this is the case if
$\eta$ is not a point of finite order on $\cE$, that is,
$N\eta\not\in\Gamma$ for any $N\in\NN$.

The space $\Theta_{n}(\Gamma)$ of the generators of the algebra
$Q_n(\cE,\eta)$ is endowed with an action of a finite group
$\wt{\Gamma_n}$ which is a central extension of the group
$\Gamma/n\Gamma$ of points of order $n$ on the curve $\cE$. It
immediately follows from the formula for the product $*$ that the
corresponding transformations of the space
$F_\alpha=S^\alpha\Theta_{n}(\Gamma)$ are automorphisms of the
algebra $Q_n(\cE,\eta)$.

\subsubsection{Bosonization of the algebra $Q_n(\cE,\eta)$}

The main approach to obtain representations of the algebra
$Q_n(\cE,\eta)$ is to construct homomorphisms from this algebra to
other algebras with simple structure (close to Weil algebras)
which have a natural set of representations. These homomorphisms
are referred to as bosonizations, by analogy with the known
constructions of quantum field theory.

Let $B_{p,n}(\eta)$ be a $\ZZ^p$-graded algebra whose space of
degree $(\alpha_1,\dots,\alpha_p)$ is of the form
$\{f(u_1,\dots,u_p)e_1^{\alpha_1}\dots e_p^{\alpha_p}\}$, where
$f$ ranges over the meromorphic functions of $p$ variables and
$e_1,\dots, e_p$ are elements of the algebra $B_{p,n}(\eta)$. Let
$B_{p,n}(\eta)$ be generated by the space of meromorphic functions
$f(u_1,\dots,u_p)$ and by the elements $e_1,\dots, e_p$ with the
defining relations
\begin{equation}
\begin{gathered}
 e_\alpha f(u_1,\dots,u_p)=
f(u_1-2\eta,\dots,u_\alpha+(n-2)\eta,\dots,u_p-2\eta) e_\alpha,\\
 e_\alpha e_\beta= e_\beta e_\alpha,\quad
f(u_1,\dots,u_p)g(u_1,\dots,u_p)=g(u_1,\dots,u_p)f(u_1,\dots,u_p)
\end{gathered}
\end{equation}

We note that the subalgebra of $B_{p,n}(\eta)$ consisting of the
elements of degree $(0,\dots,0)$ is the commutative algebra of all
meromorphic functions of $p$ variables with the ordinary
multiplication.

\begin{prop}
Let $\eta\in\cE$ be a point of infinite order. For any $p\in\NN$
there is a homomorphism $\phi_p\colon  Q_n(\cE,\eta)\to
B_{p,n}(\eta)$ that acts on the generators of the algebra
$Q_n(\cE,\eta)$ by the formula \emph:
\begin{equation}
\phi_p(f)=\sum_{1\leq\alpha\leq
p}\frac{f(u_\alpha)}{\theta(u_\alpha-u_1)\dots\theta(u_\alpha-u_p)}
e_\alpha.
\end{equation}
Here $f\in\Theta_{n}(\Gamma)$ is a generator of $Q_n(\cE,\eta)$
and the product in the denominator is of the form
$\prod_{i\neq\alpha}\theta(u_\alpha-u_i)$.
\end{prop}
\medskip
\subsubsection{Symplectic leaves}

We recall that $Q_n(\cE,0)$ is the polynomial ring
$S^*\Theta_{n}(\Gamma)$. For a fixed elliptic curve
$\cE=\CC/\Gamma$ we obtain the family of algebras $Q_n(\cE,\eta)$,
which is a flat deformation of the polynomial ring. We denote the
corresponding Poisson algebra by $q_n(\cE)$. We obtain a family of
Poisson algebras, depending on $\cE$, that is, on the modular
parameter $\tau$. Let us study the symplectic leaves of this
algebra. To this end, we note that, when passing to the limit as
$\eta\to0$, the homomorphism $\phi_p$ of associative algebras
gives a homomorphism of Poisson algebras. Namely, let us denote by
$b_{p,n}$ the Poisson algebra formed by the elements
$\sum_{\alpha_1,\dots,\alpha_p\ge0}
f_{\alpha_1,\ldots,\alpha_p}(u_1,\ldots,u_p) e_1^{\alpha_1}\ldots
e_p^{\alpha_p}$, where $f_{\alpha_1,\ldots,\alpha_p}$ are
meromorphic functions and the Poisson bracket is $$
\{u_\alpha,u_\beta\}=\{e_\alpha, e_\beta\}=0;\quad
\{e_\alpha,u_\beta\}=-2e_\alpha;\quad
\{e_\alpha,u_\alpha\}=(n-2)e_\alpha, $$
 where $\alpha\neq\beta$.

The following assertion results from Proposition 6 in the limit as
$\eta\to0$.

\begin{prop}
There is a Poisson algebra homomorphism $\psi_p\colon q_n(\cE)\to
b_{p,n}$ given by the following formula: if
$f\in\Theta_n(\Gamma)$, then 
$$
\psi_p(f)=\sum\limits_{1\leq \alpha\le
p}\frac{f(u_\alpha)}{\theta(u_\alpha-u_1)\dots\theta(u_\alpha-u_p)}
e_\alpha
$$.
\end{prop}

Let $\{\theta_i(u);i\in\ZZ/n\ZZ\}$ be a basis of the space
$\Theta_{n}(\Gamma)$ and let $\{x_i;i\in\ZZ/n\ZZ\}$ be the
corresponding basis in the space of elements of degree one in the
algebra $Q_n(\cE,\eta)$ (this space is isomorphic to
$\Theta_{n}(\Gamma)$). For an elliptic curve
$\cE\subset\PP^{n-1}$ embedded by means of theta functions of
order $n$ (this is the set of points with the coordinates
$(\theta_0(z):\ldots:\theta_{n-1}(z))$) we denote by $C_p\cE$ the
variety of $p$-chords, that is, the union of projective spaces of
dimension $p-1$ passing through $p$ points of $\cE$. Let
$K(C_p\cE)$ be the corresponding homogeneous manifold in $\CC^n$.
It is clear that $K(C_p\cE)$ consists of the points with the
coordinates 
$$
x_i=\sum\limits_{1\leq \alpha\leq
p}\frac{\theta_i(u_\alpha)}{\theta(u_\alpha-u_1)\dots
\theta(u_\alpha-u_p)} e_\alpha
$$, 
where $u_\alpha, e_\alpha\in\CC$.

Let $2p<n$. Then one can show that $\dim K(C_p\cE)=2p$ and
$K(C_{p-1}\cE)$ is the manifold of singularities of $K(C_p\cE)$.
It follows from Proposition 7 and from the fact that the Poisson
bracket is non-degenerate on $b_{p,n}$ for $2p<n$ and
$e_\alpha\neq 0$ that the non-singular part of the manifold
$K(C_p\cE)$ is a $2p$- dimensional symplectic leaf of the Poisson
algebra $q_n(\cE)$.

Let $n$ be odd. One can show that the equation defining the
manifold $K(C_{\frac{n-1}2}\cE)$ is of the form $C=0$, where $C$
is a homogeneous polynomial of degree $n$ in the variables $x_i$.
This polynomial is a central function of the algebra $q_n(\cE)$.

Let $n$ be even. The manifold $K(C_{\frac{n-2}2}\cE)$ is defined
by equations $C_1=0$ and $C_2=0$, where $\deg C_1=\deg C_2=n/2$.
The polynomials $C_1$ and $C_2$ are central in the algebra
$q_n(\cE)$.

\section{Integrable systems}

There are (at least two) different ways how to construct some
commuting elements and IS associated with the elliptic algebras.
We will start with the general statements about the commuting
elements arising from the ideas and constructions of the section 1
.

\subsection{Commuting elements in the algebras $Q_{n}(\cE;\eta)$}

Let us consider the following Weyl-like algebra $\cV_{n}$ with the
set of generators $f_1,\ldots,f_n, z_1,\ldots,z_n$ subjected to
the relations

$$ 0 =[f_i, f_j]= [z_i, z_j] = [f_i, z_j]\ (i\neq j), \ f_i z_i =
(z_i - n\eta)f_i . $$

 We assume the following commutation relations between the functions
in the variables $z_i$ and the elements $f_j$:

$$
f_jF(z_1,\ldots,z_n) = F(z_1,\ldots,z_j -n\eta,\ldots,z_n)f_j.
$$

 We should precise that the algebra $\cV_n$ is spanned as a vector space by the
elements of the form $F(z_1,\ldots,z_n)f_1^{m_1}\ldots f_n^{m_n},$
where $F$ is a meromorphic function in $n$ variables.

\begin{remark}

We should observe also that the algebra $\cV_n$ looks different
from the above-mentioned Weyl-like algebras $B_{p,n}$ but it is
isomorphic to the algebra $B_{n,n}$ and may be reduced to it by a
change of the generators. We will return below to a geometric
interpretation of the algebra $\cV_n$.

\end{remark}

 Now we take the following determinant

$$
D_0 =
\left|
\begin{array}{ccccc}
\theta_{0}(z_1)&\theta_{1}(z_1)& \ldots& \theta_{n-2}(z_1)& \theta_{n-1}(z_1)\\
\theta_{0}(z_2)&\theta_{1}(z_2)& \ldots& \theta_{n-2}(z_2)& \theta_{n-1}(z_2)\\
\vdots & \vdots& \vdots &   \vdots &         \vdots \\
\theta_{0}(z_n)&\theta_{1}(z_n)& \ldots& \theta_{n-2}(z_n)& \theta_{n-1}(z_n)
\end{array}
\right| =
$$
$$
 c \exp(z_2 + 2z_3 + \ldots +(n-1)z_n)\prod\limits_{1\leq i < j\leq n}\theta(z_i -z_j)\theta(\sum\limits_{i=1}^{n}z_i),
$$

where the constant $c$ is irrelevant for us because it will be cancelled in future computations.

Then we define the partial determinants $D_i$ replacing the $i-$th
column by the column of $f_i$,
 putting them on the place of the $n$-th column:

$$
D_i =
\left|
\begin{array}{cccccc}
\theta_{0}(z_1)&\theta_{1}(z_1)& \ldots
&^{i}&\theta_{n-1}(z_1)&f_1\\
\theta_{0}(z_2)&\theta_{1}(z_2)&\ldots& |& \theta_{n-1}(z_2)&
f_2\\ \vdots & \vdots& \vdots & \vdots & \vdots \\
\theta_{0}(z_n)&\theta_{1}(z_n)& \ldots&  |& \theta_{n-1}(z_n)&
f_n
\end{array}
\right| =
$$
$$
\sum\limits_{1\leq \al \leq n} (-1)^{\al + n}\left|
\begin{array}{ccccc}
\theta_{0}(z_1)&\theta_{1}(z_1)& \ldots&  \theta_{n-1}(z_1)\\
\theta_{0}(z_2)&\theta_{1}(z_2)& \ldots&  \theta_{n-1}(z_2)\\
\vdots & \vdots& \vdots &   \vdots &      \vdots \\
\theta_{0}(z_n)&\theta_{1}(z_n)& \ldots&  \theta_{n-1}(z_n)
\end{array}
\right|_{\al,i}f_{\al}.
$$

Here the subscription $|  \ |_{\al,i}$ means that we had to omit the 
$i$-th column and the $\alpha$- th row.

The immediate corollary of the (\ref{thm:comm}) is the following

\begin{prop}
The determinant ratios are formed a commutative family: $$ [
D_{0}^{-1}D_i, D_{0} ^{-1}D_j ] = 0. $$
\end{prop}

The result of the proposition can be expressed in an elegant way
in the form of a commutation relation of generating functions.

Let us define a generating function $T(u)$ of a variable $u \in
\CC$: $$ T(u) =  D_{0}^{-1}\sum\limits_{1\leq j \leq
n}(-1)^{j}\theta_{j}(u)D_j. $$

Then we can express the function $T(u)$, using the formulas for
the determinants of the theta-functions as

\begin{equation}\label{det}
T(u) =  D_{0}^{-1}\sum\limits_{1\leq \alpha \leq
n}(-1)^{\al}\left|
\begin{array}{ccccc}
\theta_{0}(z_1)&\theta_{1}(z_1)& \ldots& \theta_{n-1}(z_1)\\
\vdots & \vdots& \vdots &   \vdots \\
\theta_{0}(u)&\theta_{1}(u)& \ldots& \theta_{n-1}(u)\\
\vdots & \vdots& \vdots &   \vdots \\
\theta_{0}(z_n)&\theta_{1}(z_n)& \ldots& \theta_{n-1}(z_n)
\end{array}
\right|f_{\al} =
\end{equation}
$$
=
\sum\limits_{1\leq \al \leq n} {\theta(u+\sum\limits_{\beta \neq \alpha} z_{\beta})
\prod\limits_{1\leq \beta \neq
\alpha \leq n}\theta(u-z_{\beta})\over{\prod\limits_{\beta \neq
\alpha}\theta(z_{\alpha}- z_{\beta})}} {\tilde f}_{\alpha},
$$

where we denote by ${\tilde f}_{\alpha}$ the normalization
$$
{\tilde f}_{\alpha} =
{f_{\alpha}\over{\theta(\sum\limits_{i=1}^{n} z_i)}}.
$$

We remark that the commutation relations between the variables
$z_1,\ldots,z_n,{\tilde f}_1,\ldots,{\tilde f}_n$ are the same as
they were for the variables $z_1,\ldots,z_n,f_1,\ldots,f_n$.

 In this notations the proposition now reads:

\begin{prop}\label{transf}
The ``transfer-like'' operators $T(u)$ commute for different
values of the parameter $u$:
$$ [ T(u), T(v)] = 0.
$$
\end{prop}

\medskip

Now we will apply this result to a construction of the commuting
elements in the algebra $Q_{n}(\cE;\eta)$ for even $n$.

Let $n=2m$. It is known in this case that the centre
$Z(Q_{2m}(\cE;\eta))$ for $\eta$ of infinite order is generated by
Casimir elements from $S^m({\Theta}_{2m}(\Gamma))$ such that
$f(z_1,\ldots,z_m)=0$ for $z_2 = z_1 + 2m\eta$. A straightforward
computation shows that the space of such elements is
two-dimensional and has as a basis the elements of the form

$$
C_{\al} = \theta_{\al}(z_1 + \ldots z_m){\prod\limits_{1\leq i\neq
    j\leq m}}\theta(z_i - z_j - 2m\eta),
$$

where $\theta_{\al}\in {\Theta}_2(\Gamma),
\al\in{\ZZ}/{2\ZZ}.$

Let us fix an element $\Psi(z)\in \Theta_{m+5}(\Gamma)$ and two
complex numbers $a,b \in \CC$. Consider a family
$f(u)(z_1,\ldots,z_n)$ of elements from
$S^m({\Theta}_{2m}(\Gamma))$ such that $$f(u)(z_1,z_1 + 2m\eta,z_2
\ldots,z_{m-1})
$$
$$
 = \Psi(z_1 + 4m^2\eta -{{1}\over{m+5}}(a+(m-2)b + 2m(m-2)\eta))\theta(z_1 + \ldots +
z_{m-1}+a)\theta(z_1 + z_2 + b)\ldots $$ $$ \theta(z_1 + z_{m-1}+
b)\theta(z_2 - z_1 - 4m\eta)\ldots \theta(z_{m-1} - z_1 -
4m\eta)\theta(z_2 - z_1 +2m\eta)\ldots\theta(z_{m-1} - z_1 +
2m\eta)\times
$$
$$
\theta(u + z_2 +\ldots + z_{m-1}- a - b +
2m\eta)\theta(u - z_2)\ldots\theta(u - z_{m-1})\times
$$
$$
\exp(2\pi i(2(m-2)z_1 +
z_2 +\ldots + z_{m-1})){\prod\limits_{2\leq i\neq j\leq
    m-1}}\theta(z_i - z_j - 2m\eta).
 $$

The elements $f(u)$ are defined up to a linear combination of the
Casimirs $C_1, C_2$ because of the annihilation of the Casimirs on
the "diagonal" $z_1 = z_2 + 2m\eta$ ( and we are defying the
elements $f(u)$ namely on this "diagonal"!)

\begin{thm}\label{ellcomm}
In the elliptic algebra $Q_{2m}(\cE,\eta)$ the following relation
holds $$ [f(u),f(v)] = f(u)*f(v) - f(v)*f(u) = 0$$
\end{thm}

\emph{Proof} We will use the homomorphism $Q_{2m}(\cE,\eta)
\rightarrow B_{m-1, 2m}$ from the subsection 2.2.3. The element
$f(u)$ is transposed by this homomorphism into the element 
$$
\sum\limits_{1\leq\al\leq m -1}
{\theta(u + \sum\limits_{\beta\neq\al}z_\beta){{\prod\limits_{1\leq\beta\neq\al \leq m-1}
\theta(u - z_{\beta})}}\over{\prod\limits_{\beta\neq\al}\theta(z_\beta-z_\al )}}f_{\al},
$$
where we had denote by the $f_\al$ the following expression

$$ f_\al = \Psi(z_\al + 4m^2\eta -{{1}\over{m+5}}(a + (m-2)b +
2m(m-2)\eta))\times $$ $$ \theta(\sum_{\beta = 1}^{m-1}z_\beta +
a){\prod\limits_{\beta\neq\al}}\theta(z_\al + z_\beta + b)\times
 $$
 $$
\exp(2\pi i (2(m-2)z_\al + \sum_{\beta\neq\al}z_\beta))e_\al
e_1\ldots e_{m-1}.
  $$

Then, the formula (\ref{det}) and the proposition \ref{transf}
give us immediately that the images of $f(u)$ under the
homomorphism $\phi_{m-1}$ are commuting in $B_{m-1,2m}(\eta)$. It
is known that the image of the algebra $Q_{2m}(\cE,\eta)$ in
$B_{m-1,2m}(\eta)$ is the quotient of $Q_{2m}(\cE,\eta)$ over the
centre $Z = \langle C_1,C_2\rangle$. Hence the commutator $$
[f(u),f(v)] = f(u)*f(v) -f(v)*f(u) $$ belongs to the ideal
generated by the centre $Z$. To show that $[f(u),f(v)] = 0$ we can
consider the injective homomorphism into $B_{m,2m}(\eta)$ and it
is sufficient to verify that the coefficient before
$(e_1)^2\ldots(e_m)^2$ equals to zero, or ( which is equivalent)
that $$ [f(u),f(v)](z_1,z_1 + 2m\eta,z_2,z_2 +
2m\eta,\ldots,z_m,z_m + 2m\eta)=0, $$ which is achieved easily by
the direct verification.
\medskip
\begin{remark}
A family of compatible quadratic Poisson structures such that its
common linear combination is isomorphic to the Poisson structure
of the classical elliptic algebra $q_n(\cE)$ was constructed in
\cite{Od2}. The Lenard scheme enables us with a family of Poisson
commuting elements ("hamiltonians" in involution) in the algebra
$q_n(\cE)$. These elements have the degree $n$ if $n$ is impair
and $n/2$ otherwise.

 Let $n=2m$. We conjecture that the constructed in \ref{ellcomm} commuting elements in
 $Q_{2m}(\cE)$ are the quantum analogs of the commuting
 "hamiltonians" in $q_{2m}(\cE)$ generated by the Lenard scheme
 applied to the classical Casimirs ${C_0}^{(m)}, {C_1}^{(m)}$ of \cite{Od2}.
 This conjecture is true in  first non-trivial case $n=6$ when it
 is easily verified by direct computations.
\end{remark}

\subsection {Commuting elements in "bosonization" of $Q_{n,n-1}(\cE,\eta)$}

The elliptic algebra $Q_{n,n-1}(\cE,\eta)$ is commutative and the
homomorphism of "bosonization" provides a construction of a large
class of the commuting families in the corresponding Weyl-like
algebras. Some special choice of the relations between the numeric
parameters (numbers of the generators) gives us in its turn the
example of some new (to our knowledge) integrable system.

\subsubsection{Bosonization of the algebra $Q_{n,n-1}(\cE).$} The
homomorphism $\phi_p : Q_{n}(\cE,\eta)\rightarrow B_{n,p}(\eta)$
from 2.2.3 may be extended and generalized to the case of the
elliptic algebras $Q_{n,k}(\cE,\eta)$ (\cite{FO2}). The structure
of the Weyl-like algebra similar to $B_{n,p}(\eta)$ turns out to
be more complicated for $k>1$. We will use this extension in the
special case $k=n-1$ i.e. for the commutative algebra
$Q_{n,n-1}(\cE,\eta).$

Let us consider the expansion of $n\over{n-1}$ in continued
fraction to obtain the numeric parameters of the dynamical
Weyl-like algebra similar to $B_{n,p}(\eta)$:

$$
{{n}\over{n-1}} = 2 - {{1}\over{2 -\ldots -{{1}\over{2}}}}.
$$

Let ${\tilde B_{n,p_1,p_2,\ldots,p_{n-1}}(\eta)}$ be an
associative algebra generated by the following generators

$$
z_{1,1},\ldots,z_{p_1,1};e_{1,1},\ldots,e_{p_1,1};
$$
$$
z_{1,2},\ldots,z_{p_2,2};e_{1,2},\ldots,e_{p_2,2}; \vdots
$$
$$
z_{1,n-1},\ldots,z_{p_{n-1},n-1};e_{1,n-1},\ldots,e_{p_{n-1},n-1};
$$

and the additional generators
$$
t_{1,2},t_{2,3}\ldots,t_{n-2,n-1};
$$
$$
f_{1,2},f_{2,3}\ldots,f_{n-2,n-1};
$$

We will impose the following commutation relations between them:
$$ e_{\alpha,\gamma}z_{\beta,\gamma} = (z_{\beta,\gamma}
-n\eta)e_{\alpha,\gamma}; \alpha \neq \beta $$ $$ f_{\alpha,\alpha
+1}t_{\alpha,\alpha +1}= (t_{\alpha,\alpha + 1}
-n\eta)f_{\alpha,\alpha +1}. $$

All other commutators we suppose to be equal to zero.

The commutative elliptic algebra $Q_{n,n-1}(\cE,\eta)$ admits the
following "bosonization" homomorphism $\phi_{p_1,\ldots,p_{n-1}}$
which is a partial case of the general bosonization procedure for
$Q_{n,k}(\cE,\eta).$

More precisely,

$$
 \phi_{p_1,\ldots,p_{n-1}} : Q_{n,n-1}(\cE,\eta)\rightarrow
{\tilde B_{p_1,\ldots,p_{n-1}}(\eta)}
$$

can be squeezed through a homomorphism of the algebra $Q_{n-1,n}$
into an algebra of "exchange relations" with generators
$e_{\al_1,\ldots,\al_{n-1}}$(see \cite{Od1}). This "exchange"
algebra, in its turn, can be imbed in the algebra ${\tilde
B_{p_1,\ldots,p_{n-1}}(\eta)}$ with additional generators. This
statement was established in the paper \cite{FO3}. The explicit
composition map  $\phi_{p_1,\ldots,p_{n-1}}$ can be expressed like
the following ``transfer-function'' $$ \tilde T(u)= $$
\begin{equation}
{\tiny \sum\limits_{\tiny
\begin{array}{c}
1\leq \alpha_1 \leq p_1\\
\ldots\\
1\leq \alpha_{n-1} \leq p_{n-1}
\end{array}
}}
{{\theta(u-z_{\alpha_1,1})\theta(u+z_{\alpha_1,1}-z_{\alpha_2,2})\theta(u+z_{\alpha_2,2}-z_{\alpha_3,3})\ldots
\theta(u+z_{\alpha_{n-1},n-1})} \over{\prod\limits_{\tiny
\begin{array}{c}
\beta\neq\alpha_{\gamma}\\
1\leq\beta\leq p_{\gamma}\\
1\leq\gamma\leq n-1
\end{array}
} \theta(z_{\alpha_{\gamma},\gamma} -z_{\beta,\gamma})}}
\end{equation}
$$
\theta(z_{\alpha_1,1}+z_{\alpha_2,2} -t_{1,2})\theta(z_{\alpha_2,2}+z_{\alpha_3,3} -t_{2,3})\ldots
\theta(z_{\alpha_{n-2,n-2}}+z_{\alpha_{n-1,n-1}}-t_{n-2,n-1})
$$
$$
e_{\alpha_1,1}\ldots e_{\alpha_{n-1,n-1}}f_{1,2}\ldots f_{n-2, n-1}.
$$

Then the following proposition is resulted from the definition of
the combined homomorphism of the bosonization:
\begin{prop}\label{transfgen}
$$ [{\tilde T(u)},{\tilde T(v)}] = 0. $$
\end{prop}

We should observe that in the case of $p_1 = \ldots = p_{n-1} = 2$
these commuting family gives an example of  integrable system.

\medskip

\section{Some examples of elliptic integrable systems}

\subsection{Low-dimensional example: Algebra $q_{2}(\cE,\eta)$}

The elliptic integrable system arising in the \ref{ellcomm} 
becomes very transparent in
the case of $m=1$ then the commutative elliptic algebra
$Q_{2}(\cE,\eta)$ has functional dimension 2 and its Poisson
counterpart admits the Poisson morphism $$ \psi_2 : q_{2}(\cE)
\to b_2, $$ where the algebra $b_2$ are formed by the elements $$
\sum\limits_{\alpha,\beta}f_{\alpha,\beta}(z_1,z_2)e_1^{\alpha}e_2^{\beta},
$$

where $f_{\alpha,\beta}(z_1,z_2)$ are meromorphic functions and the Poisson structure in $b_2$ is given by
\begin{equation}
\{z_i,z_j \} = \{e_i,e_j \} = \{e_i,z_i \}= 0 \, \{e_i,z_j \} = - 2e_i (i \neq j), i,j=1,2.
\end{equation}
The explicit formula for this mapping (for a given theta-function
of the order 2 $f\in \Theta_2(\Gamma)$) is the following:
\begin{equation}\label{2d}
\psi_2(f) = {{f(z_1)}\over{\theta(z_1 - z_2)}}e_1 +  {{f(z_2)}\over{\theta(z_2 - z_1)}}e_2.
\end{equation}
Now, taking the basic theta-functions of the order 2 $\theta_1,
\theta_2$, let us compute the Poisson brackets between their
images $\psi_2 (\theta_1)$ and  $\psi_2 (\theta_2)$:

\begin{prop}
These theta-functions commute in $b_2$
$$
\{\psi_2 (\theta_1), \psi_2 (\theta_2) \} = \psi_2(\{ \theta_1, \theta_2 \}) = 0.
$$
\end{prop}

\subsection{SOS eight-vertex model of Date-Miwa-Jimbo-Okado}

We remind in a neccessairy form some of the SOS eight-vertex
model's ingredients (see \cite{DJMO} ) and the construction of its 
transfer-operators 
to establish the relation between them and our operators $T(u)$.

This is a statistical mechanical model which is an {\it
  Interaction-round-a- face -model} (or IRF-model )-version
of the Baxter model related to the elliptic quantum group of Felder
$E_{\tau,\eta}(sl_2)$ which was studied by the Sklyanin's
Separartion of Variables methods (under the antiperiodic boundary
conditions) in \cite{FSch} (see also \cite{FV} for the representation theory 
of $E_{\tau,\eta}(sl_2)$). We will use the results of  \cite{FSch} in
the form which we need.

The antiperiodic boundary conditions of the model are fixed by the
family of transfer-matrix $T_{SOS}(u,\lambda)$ where $u\in{\CC}$ is a
parameter and the family $T_{SOS}(u,\lambda)$ is expressed as a
(twisted) traces of the $L-$operators of the model $L_{SOS}(u,\lambda)$
defined over so-called ``auxiliary'' space of the elliptic quantum
group $E_{\tau,\eta}(sl_2)$. The $L-$operator is built out in the
tensor product of the fundamental representations of the elliptic
quantum group and is twisted by the matrix 
$\left(\begin{array}{cc}0& 1\\1& 0\end{array}\right).$

The $L-$ operator is usually represented in the form of
$2\times2-$-matrix of the form

$$
L_{SOS}(u,\lambda)=
\left(
\begin{array}{cc}
a(u,\lambda)& b(u,\lambda)\\
c(u,\lambda)& d(u,\lambda)
\end{array}
\right),
$$
where the matrix entries are meromorphic on $u$ and $\lambda$
and obey the dynamical $RLL-$commutation relations for the Felder
elliptic $R-$matrix
$$
R(u,\lambda)=
\left(
\begin{array}{cccc}
\theta(u + 2\eta)& 0& 0& 0\\
0& {{\theta(u)\theta(\lambda +2\eta)}\over{\theta(\lambda)}}&
  {{\theta(u-\lambda)\theta(2\eta)}\over{\theta(\lambda)}}& 0\\
0& {{\theta(\lambda + u)\theta(2\eta)}\over{\theta(\lambda)}}&
  {{\theta(u)\theta(\lambda - 2\eta)}\over{\theta(\lambda)}}& 0\\
0& 0& 0& \theta(u + 2\theta)
\end{array}
\right)
$$
(see {\it op.cit}).

We will have deal with the functional representations of
$E_{\tau,\eta}(sl_2)$ which are the pairs $(F,L)$ where $F$ is a
complex vector space of meromorphic functions
$f(z_1,\ldots,z_n,\lambda)$ ( or a subspace of functions which are
holomorphic in a part of the variables) and $L$ is the $L-$operator as
above. The entries of the $L-$operator are acting as the difference
operators in the tensor product $V\otimes W$ where $dim V = 2$ such
that the Felder $R-$matrix belongs to $\End(V\otimes V)$ and $W$ is an
appropriate subspace in the functional space $F$.

The Bethe ansatz method works for case of SOS model with periodic
boundary conditions after Felder and Varchenko (\cite{FV}). In the
antiperiodic case the family of the transfer matrices
$$
T_{SOS}(u,\lambda) = tr K L_{SOS}(u,\lambda) ,\ u\in {\CC}
$$
is commutative for the different parameter values:
$[T_{SOS}(u,\lambda),T_{SOS}(v,\lambda)] = 0,$
as it follows from the $RLL-$ relations by tediuos computations in
\cite{Sch}.

In the other hand it is possible to establish the explicit one-to-one correspondence
between the families of antiperiodic SOS transfer-matrices and the
auxiliary transfer-matrices $T_{aux}(u,\lambda)$ (see 4.4.3 in
\cite{Sch}). This isomorphism is established by the version of the
Separation of Variables.

The explicite expression of the auxiliar transfer-matrix is
$$
T_{aux}(u,\lambda) = \sum_{\al=1}^n{{\theta(u+z_{\al}-\lambda)}\over{\theta(\lambda)}}
\prod\limits_{1\leq \beta\neq\al\leq
  n}{{\theta(u+z_{\beta})}\over{\theta(z_{\beta} -
    z_{\al})}}(\theta(z_{\al} + \eta)T_{z_{\al}}^{-2\eta} +
\theta(z_{\al} - \eta)T_{z_{\al}}^{2\eta}),
$$
where ( to compare with our ``transfer''-operators in the elliptic
intgerable systems) we had  put in the formulas of (\cite{Sch},
ch.4.6)
$$
x_{\al} =0 ,\  \Lambda_{\al}=1 ,\  \al = 1,...,n.
$$
(The choice $\Lambda_\al = 1$ is corresponded to the case 
Separated Variables (prop. 4.36 in \cite{Sch})).

The operators $T_{z_{\al}}^{\pm 2\eta}$ are acting in the
way similar to the generators $f_\al$ above:
$$
T_{z_{\al}}^{\pm 2\eta}f(z_1,...,z_n) = f(z_1,...,z_\al \pm 2\eta,...,z_n)T_{z_{\al}}^{\pm 2\eta}.
$$

Let us make the following change of the variables
$$
f_{\al}^{\pm} = \theta(z_{\al} {\mp} \eta)T_{z_\al}^{\pm 2\eta}
$$

Now the simple inspection of the formulas shows that 

$$
T_{aux}(u, \lambda) = T_{+}(u, \lambda) + T_{-}(u, \lambda),
$$

and the the commutation results of \ref{transf} can be applied by
the following
\begin{prop}\label{Ruis}
The "transfer"-operator $T(u)$ coincides (up to an inessential numeric
factor,
depending on the normalization in the theta-function definition and
rescaling of the parameter $\eta$) with
the combinantion of the twisted traces of the auxiliary $L-$ operator
of the antiperiodic
SOS-model under the following choice of the variables:
$$
\lambda = \sum_{\al=1}^n z_{\al} 
,\ x_\al = 0, \ f_{\al}^{\pm} = \theta(z_\al {\mp}
\eta)T_{z\al}^{\pm 2\eta}.
$$
\end{prop}

The proposition \ref{transf} gives a simple proof of the commutation
relations for the twisted traces of the SOS-model in this special
case. We should observe that the commutation relations

$$
[T_{+}(u) , T_{-}(v)] + [T_{-}(u), T_{+}(v)]=0
$$
are followed from the same arguments like in \ref{transf}.

In other hand, we have got an additional argument to our belief that the
``integrability'' condition from \cite{Skew} is equivalent in some sense
to the another ``integrability'' form ( in the given case to the $RLL-$relations).

We believe also that the role of the elliptic integrable systems in
the IRF models does not restrict only to some specific conditions and
we hope to return to this point in future.

\subsection{Elliptic analogs of the Beauville-Mukai systems and Fay identity}

 Let us describe a geometric meaning of our IS's. We will start
with an observation that the classical analog of the Weyl-like
algebra $\cV_n$ given by the Poisson brackets

$$
 0 =\{f_i, f_j\}= \{z_i, z_j\} = \{f_i, z_j\}\ (i\neq j), \
\{f_i, z_i\} =  - n f_i $$ can be identified with the Poisson
algebra of the functions on the symmetric product
$S^{n}(Cone(\cE))$ of a surface represented the elliptic curve cone
( more precisely, on the Hilbert scheme $(Cone(\cE))^{[n]}$ of the
length $n$ of the points on this surface.)

The $n$ commuting elements $h_i = D_0^{-1}D_i$ in the Poisson
algebra obtaining in the \ref{transf} can be interpreted as an
elliptic version of the Beauville-Mukai systems associated with
this Poisson surface.

We will develop in more details the first interesting case ($n=3$)
of the Poisson commuting conditions (the Plucker relations from
sect.2) for this systems.

The Beauville-Mukai hamiltonians have the following form in this
case:

$$ 
H_1 = {{\det[e,\theta_1 , \theta_2]}\over{\det[\theta_0,\theta_1 ,\theta_2]}}, 
H_2 = {{\det[e,\theta_0 , \theta_2]}\over{\det[\theta_0,\theta_1 , \theta_2]}},
H_3 = {{\det[e,\theta_0 , \theta_1]}\over{\det[\theta_0,\theta_1 , \theta_2]}},
$$ 
where the
vectors-columns $e,\theta_i,\ i=0,1,2$ have the following entries
$$
e=\left(
\begin{array}{c}
e_1\\ e_2\\ e_3
\end{array}
\right),
 \theta_i =\left(
\begin{array}{c}
\theta_i (z_1)\\ \theta_i (z_2)\\ \theta_i (z_3)
\end{array}
\right),\ i=0,1,2
$$

Now the integrability condition (\ref{pluck1}) may be exprimed as
a kind of 4-Riemann identity (known also as the trisecant Fay's
identity):

\begin{equation} \label{fay}
{\tilde\theta}(v){\tilde\theta}(u + \int_{B}^A \om + \int_{D}^C)\om  = {\tilde\theta}(u +
\int_{D}^A\om ){\tilde\theta}(u + \int_{B}^C \om
){{E(A,B)E(D,C)}\over{E(A,C)E(D,B)}} +
\end{equation}
$$ 
{\tilde\theta}(u + \int_{B}^A \om ){\tilde\theta}(u + \int_{D}^C
\om){{E(A,D)E(C,B)}\over{E(A,C)E(D,B)}}, 
$$

 where $\int_{X}^Y :
\cE\otimes\cE \to \Jac(\cE) = \cE$ is a Jacobi mapping $X,Y\in
\cE, \int_{X}^Y := \int_{X}^Y \om$

for a chosen holomorphic differential $\omega$ on the elliptic
curve,

$E(X,Y) = -E(Y,X)$ is a prime form: $$
{{E(A,B)E(D,C)}\over{E(A,C)E(D,B)}} =
{{{\tilde\theta}(\int_{A}^B){\tilde\theta}(\int_{D}^C)}\over{{\tilde\theta}(\int_{A}^C){\tilde\theta}(\int_{D}^B)}}
$$

and $\tilde\theta(u)$ is an odd theta-function.
(see \cite{Fay}.)

Let $a,b,c,d$ be four points corresponded to the point $A,B,C,D$
in (\ref{fay}) and we will choose $u=a+c$.

Then the (\ref{fay}) reads as
\begin{equation}\label{pluck2}
{\tilde\theta}(a+c){\tilde\theta}(a-c){\tilde\theta}(b+d){\tilde\theta}(b-d) -
{\tilde\theta}(a+b){\tilde\theta}(a-b){\tilde\theta}(c+d){\tilde\theta}(c-d) +
\end{equation}
$$ 
{\tilde\theta}(a+d){\tilde\theta}(a-d){\tilde\theta}(c+b){\tilde\theta}(c-b) = 0,$$

where we are easily recognizing a form of the commutativity
conditions (\ref{pluck1}) for the case $n=3$ and an appropriate
choice of a linear relation between 4 points $a,b,c,d$ and
$z_1,z_2,z_3,\eta/3$ (modulo an irrelevant exponential factor
entering in the relations between the theta-functions in 
different normalizations).

\begin{remark}

1. The appearance of the Riemann-Fay relations as a commutativity
or an integrability conditions in this context looks quite natural
both from the "Poissonian " as well as from the "integrable"
viewpoints. The "Poisson-Plucker" relations and their
generalizations in the context of the Poisson polynomial
structures were studied in \cite{OR}. In the other hand, the first
links between the integrability conditions ( in the form of Hirota
bilinear identities for some elliptic difference many-body-like
systems) and the Fay formulas were established in \cite{KZW}.

2. Second remark concerns the relations between the Fay trisecant
formulas on an elliptic curve and a version of a "triangle"
relations known as the "associative" Yang-Baxter equation obtained
by Polischuk (\cite{Pol1}). This result gives an additional
evidence that the commutation relations (\ref{YB1}),(\ref{YB2}) could
be interpreted as an algebraic sort of Yang-Baxter equation. An
amusing appearance of the NC determinants in both constructions
(the Cartier-Foata determinants defined above are of course a
partial case of the quasi-determinants of Gelfand-Retakh
\cite{GR}) shows that the ideas of \cite{Skew} might be useful in
"non-commutative integrability" constructions which involve
quasi-determinants, quasi-Plucker relations etc.
\end{remark}

\section{Discussion and future problems}

We have proved that the Elliptic Algebras ( under some mild
restrictions) carry the families of commuting elements which
become in some cases the examples of Integrable Systems.

Let us indicate some questions deserving future investigations.

We are going to construct an analogs of the IS on the algebras
$Q_{2m}(\cE,\eta)$ which are corresponded to the maximal
symplectic leaves of the Elliptic Algebras $Q_{2m+1}(\cE,\eta)$
using the  bi-hamiltonian elliptic
families. This analog should quantize the bi-hamiltonian families
in ${\CC}^n$ in the case of an impair $n$.

One of the main motivations to construct the systems on $\cV_n$ 
explicitly in the terms of the determinants of $\theta-$functions
of order $n$ was the tentative to find a confirmation to the
hypothetical "separated" form of the  hamiltoninans to $n-$point
Double-Elliptic system, proposed in the paper \cite{BGOR}.

The relevance of our construction to this circle of the problems
had got an additional evidence in the paper \cite{BabTal}. Recent
discussions around of the integrability in Dijkgraaf-Vafa and
Seiberg-Witten theory provide us with the "physical" insights to
support the relation between the Beauville-Mukai and Double
Elliptic IS ( see \cite{GGG}, \cite{BH}).

The future paper (in collaboration with A. Gorsky) which should
clarify the place of the IS associated with the Elliptic Algebras
and the integrability phenomena in SUSY gauge theories are in
progress.

There are some open questions about the relations of these IS with
the Beauville-Mukai Lagrangian fibrations on the Hilbert scheme
$({\PP}^2\backslash \cE)^{[n]}$ as well as with their
non-commutative counterparts proposed in \cite{NS}. We can argue
that the proposed in \cite{ER} quantization scheme for the
fraction fields can be applied to the NC surface
$({\PP}_{S}\backslash {\cE})$ and to "extend" the deformation to
the NC Hilbert scheme $({\PP}_{S}\backslash {\cE})^{[n]}$
introducing in \cite{NS}. The proposition 3.1 from \cite{Skew}
then should in principle to give an NC integrable Beauville-Mukai
system on $({\PP}_{S}\backslash {\cE})^{[n]}$ as a NC IS
associated with the Cherednik algebras of \cite{EO}. These
questions are also in the scope of our interest.

Finally, we should mention an interesting question of
generalization of the  Beauville-Mukai IS which was introduced by
Gelfand-Zakharevich on the Hilbert scheme $(X_9)^{[n]}$ of the Del
Pezzo surface $X_9$ getting by the blow-up of 9 points on
${\PP}^2.$ These 9 points are the intersection points of two cubic
plane curves. Let $\pi$ is a Poisson tensor on $X_9$
then it can be extend to the whole ${\PP}^2$ such the extension
$\tilde\pi$ has 9 zeroes in these 9 points. The polynomial degree
of the tensor $\tilde\pi$ is 3 and it follows that this polynomial
is a linear combination of the polynomials corresponding to the
initial elliptic curves. This bi-hamiltonian structure corresponds
to the case $n=3$ studied in the paper \cite{Od2}. The open
question arises - what an algebro-geometric construction behind
the bi-hamiltonian Elliptic Poisson Algbera of \cite{Od2} in the
case of arbitrary $n$?

\section*{Acknowledgements}

{\small We would like to thank A. Gorsky, B. Enriquez, B. Khesin,
A. Levin, S. Oblezin and S. Pakuliak for helpful discussions on
the subject of the paper. A part of this work was reported as a
talk of V.R. during the Strasbourg 72-nd RCP "Quantum Groups" in
Juin 2003. We are grateful to V.Turaev and to B.Enriquez for their
invitation.

We had enjoyed the hospitality of IHES during the initial stage of
this paper. The essential part of the job was done when V.R.was a
visitor at S\'ection de math\'ematiques d'Universit\'e de
G\`en\`eve (Swiss) and at Erwin Schr\o{ö}dinger Institute for
Mathematical Physics (Wien). He is grateful to the Swiss National
Science foundation and to ESI for their invitations and support.

The work of V.R. was partially supported by CNRS (delegation
2002-2003) and partly, by the grant RFBR 01-01-00549 and by
 the grant for scientific schools RFBR 00-15-96557.

A.O. was supported partially by the grant INTAS 00 - 00055 and by
the grant for scientific school 2044-20032. A.O. also had
acknowledged the hospitality of UMR 6093, LAREMA, D\'epartement de
Math\'ematics d'Angers and LPT, Universit\'e Paris-Sud (Orsay) and
a support of NATO Scientific Collaboration Program Grant.}

\end{document}